\documentclass[16pt]{article}
\usepackage{graphicx}
\usepackage{amssymb}
\usepackage{color}
\usepackage{makeidx}         
\makeindex		     

 \let\cal\mathcal

\usepackage{graphicx}
\usepackage[curve]{xypic}
\usepackage{tikz}

\newcommand\C{{\mathbb C}}
\newcommand\Z{{\mathbb Z}}
\newcommand\N{{\mathbb N}}

 \newtheorem{theorem}{Theorem}[subsection]
\newtheorem{proposition}[theorem]{Proposition}

\newtheorem{lemma}[theorem]{Lemma}
\newtheorem{definition}[theorem]{Definition}
\newtheorem{example}[theorem]{Example}
\newtheorem{remark}[theorem]{Remark}

\usepackage{graphicx}
 
\begin{document}
\title{ Quasi-ordinary  Surface germs from the topological viewpoint}
\author{Fran\c coise Michel and  Claude Weber}

\maketitle

\begin{abstract}
Let $(W,p)$ be an analytic surface germ and let  $\nu : (W',p') \to (W,p)$ be  its normalization morphism. We can associated  to  $ (W',p')$ and $ (W,p)$ well defined links $L_{W'}$ and $L_{W}$ such that  the restriction  $\nu_L $of $\nu$ to  the link $L_{W'}$ is well defined. As $(W',p')$ is a normal  surface germ,   $p'$ is an isolated singular point and $L_{W'}$ is a  three dimensional  topological manifold.  In \cite{mich1} and \cite{mich2}, one can find a detailed proof that  $L_W$ is a topological manifold if and only if  the normalization morphism  $\nu_L$ is a homeomorphism. When the link of $(W,p)$ is not a topological manifold, one can find   in \cite{mich1} and \cite{mich2} a precise description of $L_W$ and in   \cite{mich2} a detailed study of $\nu_L$.

Here we only consider quasi-ordinary  surface germ such that the associated link is a topological manifold (such germs are, in particular, always  irreducible). We will prove, only with topological arguments, that :

1) The link of a  quasi ordinary surface germ having a bijective normalization morphism is a lens space $L(n,q)$.

2) A normal quasi-ordinary surface germ is analytically isomorphic to a cyclic quotient surface germ.

 As a consequence:
If  the link of a quasi-ordinary surface germ  is a lens space $L(n,q)$ its normalization  is a cyclic quotient surface germ. 
These results are classical (one can  refer to \cite{br1} and \cite{br2}). But,  our proofs  are new and self contained.

\end{abstract}


\bigskip


 \bigskip
 \section{Introduction}

  \begin{definition} 
   A germ $(W,0)$ of  complex  surface is {\bf  quasi-ordinary} if there exists a finite morphism 
  $\phi :(W,p) \to (\C^2,0)$ which has  a  normal-crossing discriminant.

\end{definition}

 This article  is devoted to the study of quasi-ordinary  complex surface germs  such that their associated links are three dimensional topological manifolds.
 
   \begin{definition} 
If the link associated to a complex surface germ is a topological three dimensional manifold we say the  germ has  a {\bf topologically non-singular link}. 

\end{definition}

 REMARKS:
 
 1) The key point,  in the  Hirzebruch-Jung resolution of complex surface germ ( see \cite{hi}), is the resolution of quasi-ordinary normal surface germ. It is why the singularities of  quasi-ordinary normal surface germ are often called {\bf  Hirzebruch-Jung singularities }. As explained   in \cite{mich1}, Section 4, Theorem \ref{th1}  and an explicit resolution of the singularities  $X_{n,q}=\{(x,y,z)\in \C ^3  \ s.t.\  z^n-xy^{n-q}=0\}$ where $1<n$ and $0<q<n$, $q$ prime to $n$, determine  the resolution of the  Hirzebruch-Jung singularities.
 An explicit resolution of the singularities  $X_{n,q}$ is given in Section 5.3 of \cite{mich1}. The  Hirzebruch-Jung resolution of complex surface germ  (for a presentation see  \cite{po} and Section 4 in \cite{mich1} ) is often  explicitly used  (for example see Laufer  \cite{la},   L\^e -Weber  \cite{l-w} and  Maugendre-Michel  \cite{m-m}).
 
 2) There exist quasi ordinary surface germ with non-bijective normalization morphism and topologically singular link. Here we exclude such cases. They have been studied in \cite{mich2} where  the following example is detailed: \\
 Let $Y= \{(x,y,z)\in \C ^3 \ where \  z^d-xy^d=0\}$. The normalization of $(Y,0)$ is smooth i.e. $\nu : (\C^2,0) \to (Y,0)$ is given by $(u,v) \mapsto  (u^d,v,uv)$.  Here the singular locus of $(Y,0)$ is the line $\sigma =(x,0,0), x\in \C$.   In the link $L_Y=Y\cap (D\times D\times D)$,   $N(K_{\sigma}) = L_Y \cap ( \{\vert x \vert = 1\} )$ is  a  tubular neighbourhood   of $K_{\sigma} = \sigma \cap L_Y $.  Here  $N(K_{\sigma})$ is not a topological manifold, in fact it  is a    d-curling  and  $K_{\sigma}$ is its core. Moreover $\nu$ restricted to  the line $\{(v,0), v\in \C  \}$ has degree $d$.

3) Let $n>2$ and $0<q<n-1$ be two  relatively prime integers.  The following surface germs  $X_{n,q}=\{(x,y,z)\in \C ^3  \ s.t.\  z^n-xy^{n-q}=0\}$ are quasi-ordinary but  non-normal  (when  $q \neq n-1$, the origin is not an isolated singular point).  Here (see example  \ref{elens})  we show that the link $L_X$ of $(X_{n,q},0)$ is the lens space $L(n,q)$. Such quasi-ordinary surface germs are basic examples of non-normal quasi-ordinary complex surface germs with  associated links which   are three dimensional topological  manifolds.

{ \bf{conventions and notations }}  

 $*$) We take the following  notations: All along the article, we denote by $D$ the unit disc in $\C$ i.e.  $ D=\{ z\in \C ,\vert z \vert \leq 1\}$  and by $S$ the unit circle i.e.  $ S  = \{ z\in \C ,\vert z \vert =1\}.$ For  $0<\epsilon $ we take the following notation   $ D_{\epsilon }=\{ z\in \C ,\vert z \vert \leq \epsilon \}$. \\
$**$) Let $(W,p)$ be   a quasi-ordinary  surface germ and let    $\phi :(W,p) \to (\C^2,0)$ be a finite morphism with  is a normal-crossing discriminant. Without lost of generality we may suppose that the discriminant of $\phi$ is included in $\{uv=0, (u,v) \in \C^2 \}$.  We can choose an embedding of $(W,p)$ in $(\C^n,0)$ and a compact ball $(B,0) \subset (\C^n,0)$ which is a  good compact  neighbourhood of  $(W,0)$ in $(\C^n,0)$ as defined by Durfee in \cite{du}. For sufficiently small $(\epsilon , \epsilon ')$ the intersection $(B \cap (\phi ^{-1}(D_{\epsilon}\times D_{\epsilon '}))$ is also  is a  good compact  neighbourhood of  $(W,0)$ in $(\C^n,0)$.   Moreover, to avoid an excess  of  $ `` \  \epsilon \ $",  we can transform  $(u,v)\in \C^2$  in $(u/ \epsilon   ,v /  \epsilon ')$. Then  $ (W^c,0)= ( W \cap (B\cap \phi ^{-1}(D\times D)),0)$ is  again a  good compact  neighbourhood of  $(W,0)$ in $(\C^n,0)$.  

 \begin{definition} 
  The   good compact  neighbourhood   $(W^c,0)$ in $(\C^n,0)$ describe above  (in point **) is a {\bf good compact representative of the quasi-ordinary surface germ $(W,0)$}.   The boundary of $W^c$ is  {\bf the  link  $L_W$ associated to $(W,p)$}.
  
\end{definition}

 The union of the two solid tori $L^3= ((S\times D)\cup (D\times S))$ is homeomorphic to $S^3$. 
 Let $\phi _L$ be  $\phi$  restricted to $L_W$. Then  $\phi _L$  is a ramified covering of $ L^3 $ with a ramification values  included in the Hopf link $K=((S\times \{0\} )\cup (\{0\} \times S)) \subset  L^3$. \\
 As  $L_W$ is a connected  topological manifold, $ \phi ^{-1}(S\times D)$ and $\phi^{-1} (D\times S)$   are two solid tori. It implies that $L_W$ is a lens space.
If $L_W$ is homeomorphic to $S^3$, by Mumford (\cite{mu}) the normalization of $(W,0)$, is non-singular. Non normal surface germs  with topologically  non-singular links are not classified, even when  $L_W$ is homeomorphic to $S^3$ (see Le's conjecture  (presented in \cite{l-b-s}, p. 382), partial results can be found   in \cite{l-p} and in \cite{debo}). When $L_W$ is not homeomorphic to $S^3$,  the link $L_W$, associated  to a quasi-ordinary surface germ, is a lens space $L(n,q)$ with $n>1$. The pair of integers  $(n,q)$ are  defined in Section 2, in particular $S^3$ is the lens apace $L(1,0)$.
Lens spaces associated to surface germs are studied in \cite{we} and \cite{we20}. Moreover,  our Lemma \ref{lmat} is a revisited version of Section 2 in \cite{we20}.

 In Section 4, we show the following proposition:
 
  {\bf proposition} \ref{pr1}  Let    $(W_i,p_i), i=1,2,$ be two quasi-ordinary surface germs. We suppose that their   associated links   $L_{W_i}, i=1,2$, are both   homeomorphic to  $L(n,q)$ with $0\leq q<n$.  Let $(W_i ^c,p_i)$  be a compact model of  $(W_i,p_i)$. If two    finite  analytic morphisms   $\phi_i :(W_i,p_i) \to (\C^2,0),i=1,2$,  with  a normal-crossing discriminant, determines  the same degrees  $a$ and $b$ on the  two ordered irreducible components of their ramification locus,   there exist a homeomorphism $ f: (W_1^c, p_1)  \to (W_2^c, p_2)$ such that $\phi _1= \phi_2 \circ f$.  \\  Moreover, let  $\mathring{W}_i =(W_i^c \setminus L_{W_i})$ be the interior of   $W_i^c $ and   let    $ \cal{W}_i = \mathring{W}_i \setminus (\phi _i ^{-1}( \{uv=0\})$ be the complement  of the ramification locus. The homeomorphism $f$ restricted to  $\cal{W}_1$ is an analytic isomorphism  onto  $\cal{W}_2$.

 In Section 4, we introduce the cyclic quotient surface germ $C_{n,q}$ defined as the quotient of $\C^2$ by the following action of  the cyclic group $G _n=\{\sigma \in \C, \sigma ^n=1\}$: $(z_1,z_2)\sim (\sigma ^q z_1,\sigma z_2), \ \sigma \in \cal{C}_n$.\\ 
 The definition of $C_{n,q}$ induces that the  algebra   $\cal{O} _{C_{n,q}} $ of  the analytic germs defined on $C_{n,q}$  is the sub-algebra of $\C \{z_1,z_2\}$ stable by the   action of $ G_n$. In particular $\cal{O} _{C_{n,q}} $ is generated by the monomials invariant by the action of $ G_n$ which are the monomials  $z_1^{a}z_2^{b}$ such that $n$ divises   $qa+b$.  The quotient field of $\cal{O} _{C_{n,q}} $ is also invariant by the action of $ G_n$. As the factorial ring  $\C \{z_1,z_2\}$ is integrally closed, $\cal{O} _{C_{n,q}} $  is integrally closed in its quotient field.  It proves that $C_{n,q}$  is normal.
 
 let us consider the following sequence of analytic morphisms:
 
  1) Let  $\gamma$ be  the quotient of $\C ^2$ by the  action of $ G_n$ defined just above, 
  $\gamma$ restricted to $\C ^2\setminus \{(0,0)\}$ is a regular Galois cyclic covering which ramified at $(0,0)$. \\
    2) Let us consider $\nu : C_{n,q} \to X_{n,q}$  defined by $\nu (z_1,z_2)=(z_1^n,z_2^n,z_1z_2^{n-q})$ where we denote also $(z_1,z_2)$  a representative of  $(z_1,z_2)$ in  its class in $C_{n,q}$. As $C_{n,q}$ is normal, $\nu$ is the normalisation   of $ X_{n,q}$. As  $q$ is prime to $n$, $\nu$ is bijective. So, $\nu $ restricted to the link associated to   $C_{n,q} $ is a homeomorphism. \\
  3) Let  $\pi : X_{n,q} \to \C ^2$ be the projection  $(x,y,z) \mapsto  (x,y)$  restricted to $X_{n,q}$ . So, we have $(\pi \circ \nu \circ \gamma) (z_1,z_2)=(z_1^n,z_2^n)$.  The  existence of the morphisms $\pi$ and $\pi \circ \nu $ implies that $X_{n,q}$ and  $ C_{n,q} $ are quasi-ordinary surface germs. The following sequence of finite analytic morphism will play a key role in our paper:

$$\C.{z_1} \oplus \C.z_2   \stackrel {\gamma}{\longrightarrow} C_{n,q}  \stackrel {\nu}{\longrightarrow} X_{n,q} \stackrel {\pi}{\longrightarrow}  \C.u \oplus \C.v$$

The purpose of this paper is to give, only with topological arguments,   a self-contained proof   of the following theorem (an  analytic proof  seems  more  delicate (see \cite{br2})):

 {\bf theorem} \ref{th1} 
Let    $(W,p)$ be a normal  quasi-ordinary surface germ  such that its associated link is  homeomorphic  to  the lens space  $L(n,q)$ where $0<n$ and $0\leq q<n$. Then: \\
1) If $n>1$,   $(W,p)$ is analytically isomorphic  to  cyclic quotient surface germ $(C_{n,q},0)$. If $n=1$, $(W,p)$ is smooth. \\
2)  There exists a unique, up to analytic isomorphism, finite morphism  $\psi :(W,p) \to (\C^2,0)$ of rank $n$.  Such a morphism $\psi$ has degree $1$ on its ramification locus.

{\bf Comments} 

1) Theorem \ref{th1} means that the topology of the link associated to  a normal  quasi-ordinary surface germ  determines its analytic type. 
 The famous Mumford's Theorem (\cite{mu}): `` if the link associated to a normal complex surface germ is simply connected, the surface germ is smooth", implies the particular case $n=1$ (i.e.  $L_W=L(1,0)=S^3$) of Theorem \ref{th1}.

2) By  the conic structure theorem  of Milnor (\cite{mi68}),   the local topological type of a complex analytic germ   is determined by the topology of its associated  link $L_W$. The cases where the topology of the link determines  the analytic type a complex germ are very rare.

{\bf \small Mathematics Subject Classifications (2000).}   {\small 14B05, 14J17, 32S15,32S45, 32S55, 57M45}.

\bigskip

{\bf Key words. } {\small Surface Singularities, Normalization, Lens Spaces, Discriminant, Ramified Covers }

\section{  {\bf Lens spaces }}

One can find details on lens spaces and surface singularities in \cite{we}. See also \cite{po1}.

\begin{definition} A {\bf lens space $L$} is an oriented compact three-dimensional topological manifold  which can be obtained as the union of two solid tori  $T_1\cup T_2$ glued  along their boundaries. The torus $\tau =T_1\cap T_2 $ is  the  {\bf Heegaard torus} of the  given decomposition $L=T_1\cup T_2.$
   
\end{definition}  
\begin{remark} If $L$ is a lens space, there exists  an embedded torus $\tau $ in $L$ such that $L\setminus \tau $ has two connected components  which are open solid tori $\mathring{T_i} , i=1,2$. Let $T_i, i=1,2,$ be the closure  in $L$ of the  two compact solid tori  $\mathring{T_i} $.  Of course $\tau =T_1 \cap T_2$  is a Heegaard torus for $L$. In \cite{bona}, F. Bonahon shows that  a lens space has   a unique, up to isotopy,   Heegaard torus. This  implies that the decomposition $L=T_1 \cup T_2$ is unique up to isotopy, it is  `` the"  Heegaard decomposition of $L$.
\end{remark}

A lens space $L$ with a decomposition of Heegaard torus $\tau $ can be described  as follows. The solid tori $T_i, i=1,2,$ are oriented by the orientation induced by $L$. Let $\tau _i$ be the torus $\tau $ with the orientation induced by $T_i$. By definition a meridian  $m_i$ of $T_i$ is a closed oriented  circle on $\tau_i $, essential in $\tau _i$,  which is  the boundary of a disc $D_i$ embedded in  $T_i$. A meridian of a solid torus is well defined  up  to isotopy in its  boundary. A parallel $l_i$ of $T_i$  is a closed  oriented curve on $\tau_i$ such that the intersection $m_i\cap l_i=+1$. We  also write $m_i$ (resp.  $l_i$) for  the homotopy class of $m_i$ (resp. $l_i$) in the first homotopy group of $\tau_i $. The homotopy classes  of two parallels differ by a multiple of the meridian. \\ 
We choose,   on $\tau_2 $, an oriented meridian $m_2$  and a parallel $l_2$.  
 As in \cite{we}, p. 23, we orient  a meridian $m_1$ of $T_1$ such that   $m_1= nl_2 - qm_2$ with $n\in \N$ and $q\in \Z$ where $q$ is well defined modulo $n$.  As $m_1$ is a closed curve on $\tau$, $q$ is prime to $n$. Moreover, the class of  $q$ modulo n  depends on the choice of $l_2.$ So, we can choose $l_2$ such that $0\leq q<n.$ 
 
 Let $\tau$ be a boundary component of an oriented compact three-dimensional  manifold  $M.$ Let $T$ be a solid torus given with a meridian $m $ on its boundary.  If $\gamma  $ is an essential simple closed curve  in $\tau $  there is a unique way to glue  $T$  to $M$  by an orientation reversing  homeomorphism   between  the boundary of  $T$  and  $\tau $ which send $m$ to $\gamma .$  The result of such a gluing  is  unique up to orientation preserving homeomorphism and it is called  the {\bf  Dehn filling } of $M$ associated to $\gamma $. 
 
 \begin{definition}\label{defl}

 By a Dehn filling argument,   it is sufficient to know the homotopy  class  $m_1= nl_2-qm_2$ to reconstruct $L.$  By definition {\bf the lens space $L(n,q)$}\index{lens space}  is the lens space constructed with $m_1= nl_2-qm_2$. \end{definition}

\begin{remark}\label{rlens}

 If $L$ is the union of two solid tori as described above, there   are two special cases:\
 
  1) $n=0$ if and only if  $m_1= m_2$. So,  $m_1= m_2$  if and only if  $L$ is homeomorphic to $S^1\times S^2$,
 
 2) $n=1$ if and only if  $q=0$. $q=0$ if and only if $m_1= l_2$.  So, $m_1= l_2$   if and only if  $L$ is homeomorphic to $S^3.$\
 
 \end{remark}
 
As $m_2,m_1, l_2$ are simple closed curves on the torus $\tau$ such that $m_1= nl_2 - qm_2$, $q$ is prime to $n$. To be self contained and because it is used to state   lemma \ref{lmat},  we will give a proof of the following well-known lemma. 
 
  \begin{lemma}\label{sym} 
 
 Let  suppose that $1<n$ and $1<n'$.  Two lens spaces  $L(n,q)$ and $L(n',q')$ are orientation preserving homeomorphic if and only if $n'=n$ and $q'=q$ or $q'=q^{-1} modulo \ n$.
  \end{lemma}
  
   As already said in the beginning of this section, F. Bonahon  in \cite{bona} proves  that  the torus $\tau$ is well defined up to isotopy.  The only choice we have made to define  the lens space $L(n,q)$ up to homeomorphism is the order of the  solid tori  $T_1$ and $T_2$.  The only other choice is  to begin with the torus $T_1$. Let  $\tau_1$ be  oriented as the boundary of $T_1$. Let us denote by  $\cap_i$ the intersection form on $\tau_i$. So,  $\cap_1=-\cap_2$. The above chosen closed curves $m_2,m_1$ and $l_2$ satisfy $m_2\cap_2 m_1=n>0$, $m_2\cap_2  l_2=+1$ and $m_1=nl_2-qm_2$.   Let us consider $m_i'=-m_i$.  Now  $m_1' \cap_1 m_2'= n>0$ and  $m_i'\cap_1  l_i=+1$. Moreover we have  $m_1'=-nl_2-qm_2'$. 
   
   By symmetry, we must find $a$ and $q'$ such that  $m_2'=al_1-q'm_1'$. But $m_1' \cap_1 m_2'= n=a( m_1'\cap l_1)$, it gives $a=n$. 
   Moreover, we have $m_1'=-nl_2-qm_2'$. So:
   $$ (*):m_2'=nl_1-q'm_1'=nl_1-q'(-nl_2-qm_2')=nl_2+q'nl_2+q'qm_2 ' .$$
 By construction $(m_2 ', l_2)$ gives a $\Z$-basis for $\pi _1(\tau _1, \Z)$. In   $\pi _1(\tau_1, \Z /n)$, by  the above equality (*),  implies that  $q'q=1$ in $\Z/n$.

To be self-contained, we prove  in  \ref{lens}, the  following classical result:
 if the link of a  quasi-ordinary surface germ is a topological manifold, then this link is a lens space.
 
 \begin{lemma}\label{lens} 

 Let  $\phi  :(W,p) \to (\C^2,0)$ be a finite morphism defined  on  
a  surface germ $(W,p)$. We suppose  that the link  $L_W$ of  $(W,p)$ is a topological manifold.  If  the discriminant $\Delta $ of $\phi$ is included in a normal crossing germ of curve,  then the link $L_W$ of $(W,p)$ is a lens space. The link $K_{\Gamma }$  of the singular locus $\Gamma $ of $(W,0)$,  has at most two connected components. Moreover,  $K_{\Gamma }$    is a sub-link of the  two cores of the two solid tori  of a  Heegaard decomposition of $L_W$ as  a union of two solid tori.
\end{lemma}

 {\it  Proof:}   After performing  a possible analytic  isomorphism of $(\C^2,0)$, $\Delta $ is,  by hypothesis, included in the two axes i.e. $\Delta \subset \{uv=0\}$. We follow the conventions given in Section 1.

 The union of the two solid tori $L^3= ((S\times D)\cup (D\times S))$ is homeomorphic to $S^3$. 
 Let $\phi _L$ be  $\phi$  restricted to $L_W$. Then  $\phi _L$  is a ramified covering of $ L^3 $ with a ramification locus included in the Hopf link $K=((S\times \{0\} )\cup (\{0\} \times S)) \subset  L^3$. 
 
 As  $L_W$ is a  connected topological manifold,  $ \phi ^{-1}(S\times D)$ and $\phi^{-1} (D\times S)$   are two solid tori. It implies that $L_W$ is a lens space.  Moreover  $ \phi ^{-1}(S\times \{0\})$ (resp. $\phi^{-1} (\{0\}\times S)$)   is  connected, it is the core of the solid torus $ \phi ^{-1}(S\times D)$ (resp. $\phi^{-1} (D\times S)$). The singular locus   $\Gamma $ of $(W,0)$ is included in the singular locus $ \phi ^{-1}(D\times \{0\}) \cup \phi^{-1} (\{0\}\times D)$)  of $\phi$.

   {\it End of proof.}
   
   \begin{example}\label{elens} Let $n$ and $q$ be two  relatively prime strictly positive integers.  We suppose that  $0<q<n$. Let $X_{n,q}=\{(x,y,z)\in \C ^3  \ s.t.\  z^n-xy^{n-q}=0\}$. The link $L_{X_{n,q}}$ of $(X_{n,q},0)$ is the lens space $L(n,q)$. Moreover $(X_{n,q},0)$ is normal if and only if $q=n-1$.
   \end{example}

   As, $(X_{n,q},0)$ is a  hypersurface  in $\C^3$, $(X_{n,q},0)$ is normal if and only if $0$ is an isolated singular point. But,  $0$ is an isolated singular point if and only if $n-q=1$. Moreover, when  $n\neq 1+q$ the boundary of the Milnor fiber of $f (x,y,z)= z^n-xy^{n-q}$ which  is determined in \cite{m-p-w} is never homeomorphic to $L_{X_{n,q}}$.
   
    Let $\pi : (X_{n,q},0)\to (\C^2,0)$ be the projection  $(x,y,z) \mapsto  (x,y)$ restricted to $X$. 
   The discriminant  $\Delta $ of $\pi$ is equal to $ \{uv=0\}$. By Lemma  (\ref{lens}),    $L_{X_{n,q}}$ is a lens space.  \\
    Let   $ T_1=\pi ^{-1}(S\times D)$ and $ T_2=\pi^{-1} (D\times S)$  be the  two solid tori of the Heegaard decomposition of $L_{X_{n,q}}$. \\
  As $n$ and $q$ are  relatively  prime, the equation of  $X_{n,q}$ implies that   $ D_1= \pi ^{-1}(\{a\} \times D)$  and  $ D_2 = \pi ^{-1}(D\times \{b\})$ are connected discs. The oriented  boundary  $m_i,\  i=1,2, $ of $D_i$ is a meridian of $T_i$.  \\
  We choose $c\in \C$  such that $c^n=ab^{n-q}$. Let $l_2=  \{z=c\} \cap \pi ^{-1}(S \times S)$. On the torus  $ \tau = \pi ^{-1}(S \times S)$, oriented as the boundary of $T_2$, we have  $m_2\cap m_1=+n$ and we choose the orientation of $l_2$ such that  $m_2\cap l_2=+1$. Then, we also have that  $m_1\cap l_2=n-q$. \\
   We write $m_1=\alpha l_2 +\beta m_2$. Then $ n= m_2\cap m_1 = \alpha$ and  $ n-q=m_1\cap l_2=\beta $.  So $m_1= n (l_2+ m_2) -qm_2$ and   $l_2'=l_2+ m_2$ is another parallel such that $m_1=n l_2'-q m_2$. By definition $L_{X_{n,q}}$ is the lens space $L(n,q)$. \\

 \section{Finite morphisms with smooth discriminant}

The following Lemma is used in the construction of the Hirzebruch-Jung resolution of complex surface germs. It is why we write  here a self-contained proof of it.

\begin{lemma}\label{lsmooth} 

  Let  $\phi  :(W,p) \to (\C^2,0)$  be a finite morphism,   of generic degree  $n$,  defined  on  
a normal surface germ $(W,p).$   If the discriminant of $\phi $ is a smooth germ of curve, then     $(X,0)$ is analytically isomorphic to $(\C^2,0)$ and   $\phi $ is analytically  isomorphic to  the map from $(\C^2,0)$  to $(\C^2,0)$ defined by $(x,y)  \mapsto (x,y^n).$

\end{lemma}

{\it  Proof:}  After performing  an  analytic  automorphism  of $(\C^2,0)$, we can  choose coordinates such that $\Delta=\{v=0\} $. \\
 As explained in the previous subsection, the union of the two solid tori $L^3= ((S\times D)\cup (D\times S))$ is homeomorphic to $S^3$.  \\
 Let $\phi _L$ be  $\phi$  restricted to $L_W$. Then  $\phi _L$  is a ramified covering of $ L^3 $ with a ramification locus included in the trivial knot  $K=(S\times \{0\}) \subset  L^3$. 
  As $(W,p)$ is normal, $L_W$ is a  connected topological manifold. Then,   $T_1= \phi_L ^{-1}(S\times D)$ and $T_2=\phi_L^{-1} (D\times S)$   are two solid tori. It implies that $L_W$ is a lens space.  Moreover  $ \phi_L ^{-1}(S\times \{0\})$ (resp. $\phi_L^{-1} (\{0\}\times S)$)   is  connected, it is the core of the solid torus $ \phi_L ^{-1}(S\times D)$ (resp. $\phi_L^{-1} (D\times S)$). The singular locus of $\phi_L$ is included in  $ \phi_L ^{-1}(S\times \{0\}) $. Let $(a,b)\in S\times S$. As $(o,b)$ is a regular value of  $\phi_L$, $ \phi_L ^{-1}(D\times \{b\})$   is  the disjoint union of  $n$  discs, where $n$ is the generic order of $\phi$. Let $D_2$ be one of these $n$ discs and $m_2$ be the oriented  boundary of $D_2$. The order of $\phi_L $ restricted  to $m_2$ is equal to $1$. If $d$ is the number of points of $\phi_L^{-1} ( \{a\}  \times \{0\})$, then $ \phi_L ^{-1}(\{a\} \times D)$   is the  disjoint union of $d$  discs, let $D_1$ be one of them and $m_1$ be the  oriented boundary of $D_1$.\\
  As the order of $\phi_L $ restricted  to $m_2$ is equal to $1$,   the intersection $D_1 \cap D_2$ has only one point.
 As $m_2$ is the boundary of one of the $n$ connected components of  $ \phi_L ^{-1}(D\times \{b\})$,  there is an orientation for $m_1$ such that $m_2 \cap m_1 =+1$ and $m_1$ can be a parallel  $l_2$ for  $T_2$. This  is the case 2) in  \ref{rlens},  so  the link $L_W$ of $(W,p)$ is  the 3-sphere $S^3$. 
 As $(W,p)$ is normal, by Mumford \cite{mu}, $(W,p)$ is a smooth surface germ i.e  $( W,p)$ is  analytically  isomorphic to $(\C^2,0)$. The first part of the statement has been proved.

(*) Moreover $\phi_L^{-1} (S \times \{0\}) \cup  (\{ 0\}  \times S)$ is the union of the cores of $T_1$ and $T_2$. Then, {\bf  $(S \times \{0\}) \cup  (\{ 0\}  \times S)$ is a Hopf link in the 3-sphere $L_W$}.

From now on, $\phi  :(\C^2,0) \to (\C^2,0)$ is a finite morphism and its discriminant locus is $\{v=0\}$.
Let us write  $\phi=(\phi_1, \phi_2)$. The  link of the zero locus of the function germ $$(\phi_1.\phi_2)  :(\C^2,0) \to (\C ,0)$$ is the link  describe   above (see (*)), i.e.  it is   a Hopf link.  The  function $(\phi_1 .\phi_2) $ reduced is analytically  isomorphic to $(x,y)  \mapsto (xy).$   
But,  $ (\phi_L ^{-1}(\{a\} \times D)) \cap ( \phi_L ^{-1}(D\times \{b\})) =  \phi _L^{-1}( a \times b)$ has exactly $n$ points. On  the torus $\tau =T_1\cap T_2$ each of the $d$ connected components of the boundary of $ \phi_L ^{-1}(\{a\} \times D)$ cuts each of the $n$ connected components of the boundary of $ \phi_L ^{-1}(D\times \{b\}) $ in one points. So, $d=1$ because   $\phi _L^{-1}( a \times b)$ has $nd=n$ points, Then, $ \phi_L ^{-1}(\{a\} \times D)=D_1$ is a disc.  It implies that $\phi_1 $ is reduced because its  Milnor fiber  is diffeomorphic to  $D_1$. So, $\phi_1$ is isomorphic to  $x$ and $\phi_2=y^m, m\geq 1$. \\
But, $m=n$  because the Milnor fiber of $\phi_2 $ is diffeomorphic to   $  \phi_L^{-1}(D \times \{b\})$  which is the disjoint union of $n$ discs. This completes   the proof that $\phi_2$ is isomorphic to $y^n$ and  $\phi=(\phi_1, \phi_2)$ is isomorphic to $(x,y^n)$.

{\it End of proof}.

\section{  {\bf Characterization  of  coverings  with  have the Hopf link  as  ramification locus }}

 Let  $\phi  :(W,p) \to (\C^2,0)$ be a finite morphism defined  on  
a  surface germ $(W,p)$. We suppose  that the link  $L_W$ of  $(W,p)$ is a topological manifold and  the discriminant $\Delta $ of $\phi$ is a normal crossing germ of curve. Following the conventions  of \ref{lens} ,  we have $\Delta = \{uv=0\}$ (the case with a smooth discriminant has been treated in the preceding section).

 The union of the two solid tori $L^3= ((S\times D)\cup (D\times S))$ is homeomorphic to $S^3$. 
 Let $\phi _L$ be  $\phi$  restricted to $L_W$.  As  $L_W$ is a  connected topological manifold,  $T_1= \phi _L^{-1}(S\times D)$ and $ T_2= \phi _L^{-1} (D\times S)$   are two solid tori. 
 So, the ramification locus of $\phi _L$ is  exactly the  union of the two knots  $K_1= \phi _L ^{-1}(S\times \{0\}) $ and $K_2=\phi _L^{-1} (\{0\}\times S)$.    Moreover,  $K_1\cup K_2 $  are  the  two cores of the two solid tori  of a  Heegaard decomposition of  the lens space $L_W$ as  the union $T_1 \cup T_2$. As in Section 2, the torus  $\tau =T_1\cap T_2$ is  oriented as the boundary of $T_2$. \\
 Then  $\phi _L$  is a ramified covering of $ L^3 $ with  ramification  the Hopf link $K=((S\times \{0\} )\cup (\{0\} \times S)) \subset  L^3$. Let $\phi_{h} :\pi_1(\tau) \to \pi_1(S\times S)$ be the homomorphism induced by $\phi_L$. Let $e_1$ be the  first homotopy class of $\phi  _L(K_2)=(\{0\}\times S)$ and $e_2$ be the  first homotopy class of $\phi  _L(K_1)=(S \times \{0\} )$. The   meridians of the Hopf link are   unique up to isotopy.  The Hopf link can be numbered to define  the ordered basis  $(e_1,e_2)$ of $ \pi_1(S\times S)$.

\begin{lemma}\label{lmat} 
 Let  $\phi  :(W,p) \to (\C^2,0)$ be a finite morphism defined  on  a  surface germ $(W,p)$ such that  $\Delta = \{uv=0\}$ is its discriminant locus.
 Let us suppose that $L_W$ is  the  lens space $L(n,q)$ with $0\leq q<n$. Let  $\phi  _L  : L_W  \to L^3$  be the ramified covering  restriction of $\phi$ to $L_W$.  The  set of the  ramification  values of $\phi_L$ is  the Hopf link $K=((S\times \{0\} )\cup (\{0\} \times S)) \subset  L^3$.  Then  $T_1 =\phi _L ^{-1}(S \times D)$ and  $T_2 =\phi _L ^{-1}(D \times S)$ are the two solid tori of the heegaard decomposition $L_W=T_1\cup T_2$. 
  We ordered the coordinates of $\C ^2$   such that  the meridian $m_2$  and the parallel $l_2$ of $T_2$, satisfied $m_1=nl_2-qm_2$.   In the ordered basis $(m_2,l_2)$ of $\pi_1 (\tau)$   and  the ordered basis  $(e_1,e_2)$ of $ \pi_1(S\times S)$,  the matrix of $\phi_{h}$ is equal to:
  $$A= \left( \begin{array}{cc}
an & aq\\
0& b
\end{array}\right),
$$
Where $a$ is the degree of $\phi  _L $ on  the core $K_1$  of $T_1$ and $b$ the degree of $\phi  _L $ on  the core $K_2$  of $T_2$.

\end{lemma}

{\bf{Remark}}

 As in \ref{lmat}, let   $\phi  _L  : L_W  \to L^3$  be the ramified covering  restriction of $\phi$ to $L_W$.  The  set of the  ramification  values of $\phi_L$ is  the Hopf link $K=((S\times \{0\} )\cup (\{0\} \times S)) \subset  L^3$.  Let $a$ be the degree of $\phi_L$ restricted to  $K_1 =\phi _L ^{-1}(S \times \{0\})$ and $b$ be the degree of $\phi_L$ restricted to  $K_2 =\phi _L ^{-1}(\{0\} \times S)$.  Then  $T_2 =\phi _L ^{-1}(D \times S)$ is one of  the two solid tori of the heegaard decomposition $L_W=T_1\cup T_2$.  We can choose an order  for the  coordinates of $(\C^2,0)$ such that  the meridian $m_2$ and the parallel $l_2$ of $T_2$, are such that  $m_1=nl_2-qm_2.$ If $q=0$,  we have $n=1$ and $m_1=l_2$.    As proved in \ref{sym},  if  $0<q$ and if $ T_2' =\phi _L ^{-1}(D \times S)$ is the other torus of the heegaard decomposition $L_W=T_1' \cup T_2'$,  the meridian $m_2'$  and the parallel $l_2'$ of $T_2'$, are such that  $m_1'=nl_2'-q'm_2'$ where $0<q'<n$ and $qq'=1$ modulo $n$. For that choice of order, the matrix associated to $\phi_L$ is equal to $A'$ where:
  $$A'= \left( \begin{array}{cc}
an & aq'\\
0& b
\end{array}\right).
$$

 So, up to the order of  the coordinates of  $\C ^2$, the matrix $A$  associated to  $\phi_L$  is unique.  Lemma \ref{lmat} allows us to state the following definition.

  \begin{definition}\label{dmat} 
   Let  $\phi  :(W,p) \to (\C^2,0)$ be a finite morphism defined  on  a  surface germ $(W,p)$ such that  $\Delta = \{uv=0\}$ is its discriminant locus.  We suppose  that the link  $L_W= =\phi _L ^{-1}(S \times D) \cap \phi _L ^{-1}(D \times S)$ of  $(W,p)$ is the lens space $L(n,q)$ with $0<q<n$. The restriction of $\phi$ to the link associated to $(W,0)$ is the ramified covering   $\phi  _L  : L_W  \to L^3$    with  ramification  the Hopf link $K=((S\times \{0\} )\cup (\{0\} \times S)) \subset  L^3$.  The {\bf matrix associated to $\phi$} is the matrix 
  $$A= \left( \begin{array}{cc}
an & aq\\
0& b
\end{array}\right),
$$
Where $a$ is the degree of $\phi  _L $ on the core of $T_1$ and $b$ the degree of  $\phi _L$ on the core of $T_2$.

\end{definition}

{\it proof of \ref{lmat}}

 We consider the following topological invariants of $\phi _L$:
 The order $a$ of $\phi_L$ on $K_1$ and the order $b$ of $\phi_L$ on $K_2$. So, $\phi_L ^{-1}(1,0)$ has $a$ points and $\phi _L ^{-1}(\{ 1\} \times D)$ consists in $a$  oriented discs, let $D_1$ be one of these $a$ discs. For the same reason,  $\phi  _L^{-1}(0,1)$ has $b$ points and $\phi _L ^{-1}(D \times \{1\})$  consists in $b$  oriented discs, let $D_2$ be one of these $b$ discs. Let $m_i$ be the boundary  of $D_i$ where $m_2$ is oriented as the boundary of $D_2$  and $m_1$ is oriented in such a way that the intersection $m_1\cap m_2$ on the oriented torus $\tau$ is  $ n>0$. So, the degree of $\phi _L$ on $m_2$ is equal to $an$, the degree of $\phi _L$ on $m_1$ is equal to $bn$ and the  generic order of  $\phi _L$ is  equal to  $abn$.

 A parallel $l_i$ of $T_i$  is a closed  oriented curve on $\tau_i$ such that the intersection $m_i\cap l_i=+1$. We  also write $m_i$ (resp.  $l_i$) for  the homotopy  class of $m_i$ (resp. $l_i$) in  $\pi_1 (\tau)$.  As before we  choose $l_2$ such that $m_1=nl_2-qm_2$ where $0\leq q<n$. By the  choice of $0\leq q<n$, $l_2 \in  \pi_1 (\tau)$ is unique. So $\phi_1 (l_2)=qa\ e_1+b\ e_2$.
 Moreover $n$ and $0\leq q<n$ are given by the homeomorphism class of the lens space $L_W$. The degree  $a$ of $\phi _L$ on the core $K_1$ of the solid torus $T_1$  is given by the well defined matrix $A$, it is the g.c.d. of its first line. The degree  $b$ of $\phi _L$ on the core $K_2$ of the solid torus $T_2$  is  the coefficient $2\times 2$ of 
 the well-defined matrix $A$.
 
 {\it{end of proof of \ref{lmat}}}
 
 {\bf{Remark}}

  Let  $\phi  :(W,p) \to (\C^2,0)$ be a finite morphism defined  on  a  surface germ $(W,p)$ such that  $\Delta = \{uv=0\}$ is its discriminant locus.
 Let us suppose that $L_W$ is  homeomorphic to the sphere $S^3$. Let  $\phi  _L  : L_W  \to L^3$  be a ramified covering   with  ramification  the Hopf link $K=((S\times \{0\} )\cup (\{0\} \times S)) \subset  L^3$. Let $a$ be the degree  of $\phi_L$ on  $K_1= \phi _L ^{-1}(S\times \{0\}) $ and $b$ be the degree of $\phi_L$ on  $K_2=\phi _L^{-1} (\{0\}\times S)$. 
  Let $m_1$ be the boundary of one of the $a$ discs  $\phi _L ^{-1}(\{ 1\} \times D)$. Let $m_2$ be  the oriented boundary of one of the $b$ discs   $\phi  _L^{-1}(0,1)$.  We oriente $m_1$ such that $m_1\cap m_2$ is positive. But,  $L_W$ is homeomorphic to $S^3$ if and only if  $n=m_1\cap m_2=1$ and $q=0$. This remark  is just a particular case  of \ref{lmat} where  $a.b$ is the generic order of $\phi_L$ and the matrix associated to $\phi$ is equal to:

  $$A_0= \left( \begin{array}{cc}
a & 0\\
0& b
\end{array}\right),
$$

   \begin{example}\label{ematlens} Let $n$ and $q$ be two  relatively prime strictly positive integers.  We suppose that  $0<q<n$. Let $X_{n,q}=\{(x,y,z)\in \C ^3  \ s.t.\  z^n-xy^{n-q}=0\}$. 
   Let  $\pi : X_{n,q} \to \C ^2$ be the projection  $(x,y,z) \mapsto  (x,y)$  restricted to $X_{n,q}$. As computed in \ref{elens},   the link $L_{X_{n,q}}$ of $(X_{n,q},0)$ is the lens space $L(n,q)$. Here the degrees $a$ and $b$ are equal to $1$. By  Lemma \ref{lmat}, the matrix associated to $\phi$ is equal to:
   
     $$A_1= \left( \begin{array}{cc}
n & q\\
0& 1
\end{array}\right),
$$

   \end{example}

  \begin{definition}\label{dquotient}  
  
     Let $n$ and $q$ be two  relatively prime strictly positive integers.  We suppose that  $0<q<n$.  {\bf The cyclic quotient surface germ } $C_{n,q}$  is   the quotient of $\C^2$ by the following action of $\cal{C} _n=\{\sigma \in \C, \sigma ^n=1\}$ : $(z_1,z_2) \sim (\sigma ^q z_1,\sigma z_2),  \sigma \in \cal{C} _n$.
 
\end{definition}

   \begin{example}\label{equotient}
   
    Let $n$ and $q$ be two  relatively prime strictly positive integers.  We suppose that  $0<q<n$. Let us consider the cyclic quotient surface germ $C_{n,q}$. As remark in the introduction $C_{n,q}$ is normal.
   Let us consider $\nu : C_{n,q} \to X_{n,q}$  defined by $\nu (z_1,z_2)=(z_1^n,z_2^n,z_1z_2^{n-q})$ where we denote also $(z_1,z_2)$  a representative of  $(z_1,z_2)$ in  its class in $C_{n,q}$. As $C_{n,q}$ is normal, $\nu$ is the normalization   of $ X_{n,q}$. As  $q$ is prime to $n$, $\nu$ is bijective.  So, $\nu $ restricted to the link associated to   $C_{n,q} $ is a homeomorphism. As the lens space $L_(n,q)$ is the link associated to $ X_{n,q}$, the  link associated to $C_{n,q} $  is also  the lens space $L(n,q)$.

   Let  $\pi : X_{n,q} \to \C ^2$ be the projection  $(x,y,z) \mapsto  (x,y)$  restricted to $X_{n,q}$. The  discriminant of  the finite morphism  $\pi \circ \nu: C_{n,q} \to \C ^2 $ ls also  equal to $\Delta = \{uv=0\}$ and  $\pi \circ \nu$   restricted to the link of $C_{n,q}$ has the same topology than $\pi$ restricted to the link of $X_{n,q}$.  In particular  $\pi \circ \nu $ has degree $1$ on the two irreducible components of its ramification locus $(\pi \circ \nu)^{-1} (\Delta)$ i.e  the degree of   $\pi \circ \nu$  restricted  to  $ (\pi  \circ \nu)^{-1}(S\times \{0\}) $ is equal to $1$,   and  the degree of   $\pi \circ \nu$  restricted  to $(\pi \circ  \nu )^{-1} (\{0\}\times S)$ is also equal to $1$.
    By  Lemma \ref{lmat} and Example \ref{ematlens},  the matrix associated to $\pi \circ \nu$  is equal to the matrix $A_1$ associated to $\pi$:
   
     $$A_1= \left( \begin{array}{cc}
n & q\\
0& 1
\end{array}\right).
$$

   \end{example}

 \section{  {\bf Finite morphism with normal-crossing discriminant }}

\begin{proposition}\label{pr1} 
 Let    $(W_i,p_i), i=1,2,$ be two quasi-ordinary surface germs. We suppose that their   associated links   $L_{W_i}, i=1,2$, are both   homeomorphic to  $L(n,q)$ with $0\leq q<n$.  Let $(W_i ^c,p_i)$  be a compact model of  $(W_i,p_i)$. If two    finite  analytic morphisms   $\phi_i :(W_i,p_i) \to (\C^2,0),i=1,2$,  with  a normal-crossing discriminant, determines  the same degrees  $a$ and $b$ on the  two ordered irreducible components of their ramification locus,   there exist a homeomorphism $ f: (W_1^c, p_1)  \to (W_2^c, p_2)$ such that $\phi _1= \phi_2 \circ f$.  \\  Moreover, let  $\mathring{W}_i =(W_i^c \setminus L_{W_i})$ be the interior of   $W_i^c $ and   let    $ \cal{W}_i = \mathring{W}_i \setminus (\phi _i ^{-1}( \{uv=0\})$ be the complement  of the ramification locus. The homeomorphism $f$ restricted to  $\cal{W}_1$ is an analytic isomorphism  onto  $\cal{W}_2$.

\end{proposition}

{\it{ proof of \ref{pr1}}}

 By  Lemma \ref{lmat}, $\phi_1$ and $\phi_2$ have the same associated matrix $A$. By the classification of topological regular coverings,  the regular coverings $\phi_{i,\tau}$,  defined by the restriction of  $\phi _i, i=1,2,$  to the torus  $ \tau _i= \phi  _i ^{-1}(S \times S)$,  are homeomorphic coverings by a homeomorphism $f_{\tau}$. \\
   (*) The chosen  basis of $\pi_1 (\tau _i)$  (se \ref{lmat}) to write $A$ implies that $f_{\tau }$ sends the meridians of $\phi  _1 ^{-1}(D \times S)$ onto the meridians of $\phi  _2 ^{-1}(D \times S)$ (resp. the chosen  parallels of $\phi  _1 ^{-1}(D \times S)$ onto the parallels  of $\phi  _2 ^{-1}(D \times S)$) and has   degree $1$ on all these simple curves. The same is true for meridians and parallels of  the solid tori $\phi  _1 ^{-1}(S \times D)$ and $\phi  _2 ^{-1}(S \times D)$.  \\
   For all $u \in S $ the homeomorphism $f_{\tau}$ is defined on the $a$ meridians $\phi _1 ^{-1}(\{ u\} \times S)$ of the solid torus $\phi  _1 ^{-1}(S \times D)$. Let $m_u$  be one of these $a$ meridians. The definition of the  matrix $A$  fixes the behavior of $f_{\tau }$ on meridians and parallels as stated in (*) just above.  So, $f_{\tau } (m_u)$ is a  meridian of  $\phi  _2 ^{-1}(S \times D)$. Then  $f_{\tau }$ extends continually  to the $a$ points  of $\phi _{1} ^{-1}( u  \times 0)$, it is the topological conic continuous  extension  of $f_{\tau }$ to the origin of  the compact disc $\phi _1 ^{-1}(\{ u\} \times  D)$. For the same reason,  $f_{\tau}$ extends continually  to the $b$ points  of $\phi _{1} ^{-1}( 0  \times v), v \in S$. Let us denote   by $f_L: L_{W_1} \to L_{W_2}$ this extension of $f_{\tau}$ to the link $l_{W_1}$. But $f_{\tau }$  has degree  $1$ on the boundaries of all the meridian discs. So, $f_L $ is a homeomorphism of ramified coverings.    This homeomorphism $f_L$  has  a topological conic continuous  extension  $f$ to  $W_1^c$. So, we have constructed a homeomorphism $f$:
 
 $$ f: (W_1^c, p_1)  \to (W_2^c, p_2)$$
 
  such that $\phi _1= \phi_2 \circ f$. 
 
  Let  $\mathring{W}_i =(W_i^c \setminus L_{W_i})$ be the interior of   $W_i^c $.   
But $\phi_i, i=1,2$,  restricted to     $ \cal{W}_i = \mathring{W}_i \setminus (\phi _i ^{-1}( \{uv=0\})$
are  analytic  regular coverings. if $p\in \cal{W}_1 $ there exists an open neighbourhood $V_1$ of $p$ in $\cal{W}_1 $ such that $\phi _1$ restricted to $V_1$ is an analytic isomorphism. Then $\phi _2$ restricted to $V_2=f(V_1)$ is also an analytic isomorphism. We can locally  inverse  $\phi _2$ on $\phi  _1 (V_1)$. Then,  $f$  restricted to  $V_1$   is equal to    $  \phi _2^{-1} \circ \phi _1$ which is an analytic isomorphism. Let $f_0$ be $f$ restricted to $\cal{W}_1 $. We have proved that $f_0$ is an analytic isomorphism. But,  in general, $f$ is not analytic around the points of the ramification locus. 

 {\it{  end of proof of \ref{pr1}}}

\begin{theorem}\label{th1}

Let    $(W,p)$ be a normal  quasi-ordinary surface germ  such that its associated link is  homeomorphic  to  the lens space  $L(n,q)$ where $0<n$ and $0\leq q<n$. Then: \\
1) If $n>1$,   $(W,p)$ is analytically isomorphic  to  cyclic quotient surface germ $(C_{n,q},0)$. If $n=1$, $(W,p)$ is smooth. \\
2)  There exists a unique, up to analytic isomorphism, finite morphism  $\psi :(W,p) \to (\C^2,0)$ of rank $n$.  Such a morphism $\psi$ has degree $1$ on its ramification locus.

\end{theorem} 

{\bf Comments} Theorem \ref{th1} means that the topology of the link associated to  a normal  quasi-ordinary surface germ  determines its analytic type. The  case $n=1$  (i.e. when $L_W$ is homeomorphic to $S^3$) is  the  Mumford's Theorem (\cite{mu}) for normal quasi-ordinary surface germs.

 {\it{  proof of \ref{th1}}}
 
 As $(W,p)$ is quasi ordinary, there exist a finite morphism  $\phi :(W, p) \to (\C^2,0)$, such that  $\Delta = \{uv=0\}$ is   its discriminant locus. The  trivial case with  a smooth discriminant is treated in \ref{lsmooth}.
 By Lemma \ref{lmat} the matrix associated to $\phi$ is equal  to:
   $$A= \left( \begin{array}{cc}
an & aq\\ 
0& b
\end{array}\right),
$$
 Where    $a$ and $b$ are the generic degree of $\phi$ on the  two ordered irreducible components of its  ramification locus.
 
 Let $\pi \circ \nu: (C_{n,q},0)  \to (\C ^2,0) $ be the finite morphism defined in \ref{equotient}. The matrix associated to   $\pi \circ \nu$ is equal to    $$A_1= \left( \begin{array}{cc}
n & q\\
0& 1
\end{array}\right),
$$

Let $\alpha: (\C ^2,0)   \to (\C ^2,0) $ be defined by $\alpha (u,v)=(u^{a}, v^{b})$, the matrix associated to $\alpha $ is the diagonal matrix $A_0$ (see \ref{lmat} when $q=0$). So $A$ is the matrix associated to $$  \alpha \circ \pi \circ \nu: (C_{n,q},0)  \to (\C ^2,0) .$$

 In particular  $(\phi :(W, p) \to (\C^2,0))$ and $(( \alpha \circ \pi \circ \nu): (C_{n,q},0)  \to (\C ^2,0)) $ determine  the same degrees  $a$ and $b$ on the  two ordered irreducible components of their ramification locus. Proposition \ref{pr1} proves the existence of a homeomorphism $f$ between two compact models $(W^c,p)$ and $(C_{n,q}^c,0)$,   of the surface germs $(W,p)$ and $(C_{n,q},0)$.  The constructed   $ f: (W^c, p)  \to (C_{n,q}^c, 0) $ is  such that $\phi =  \alpha \circ \phi \circ \nu \circ f$.  
 
   On the other hand the interior $\mathring{W}=(W^c \setminus L_W)$ of $W^c$  (resp.  the interior $\mathring{C}_{n,q}=(C_{n,q}^c \setminus L_{C_{n,q}})$ of $C_{n,q}^c$)  is a good open analytic neighbourhood of $p$ in $(W,p)$ (resp. of $0$ in  $(C_{n,q},0)$).
   
   Let   $ \cal{W}= \mathring{W}_i \setminus (\phi ^{-1}( \{uv=0\})$,  (resp   $ \cal{C} _{n,q} = \mathring{C}_(n,q) \setminus ( ( \alpha \circ \pi \circ \nu) ^{-1}( \{uv=0\}))$) be the complement of the ramification locus.  By Proposition \ref{pr1}, the homeomorphism $f$ restricted to  $\cal{W}$ is an analytic isomorphism  onto  $ \cal{C} _{n,q}$.

 But $f$ is continuous on $W^c$, $(W,p)$ is normal  and the ramification locus of $\phi$ is a  one-dimensional  sub-analytic space of $(W,p)$, then $f$ restricted to $( \mathring{W}, p)$ is  analytic  onto $\mathring{C}_{n,q}$.  But $(C_{n,q},0)$ is also normal. Then $f^{-1}$ is also analytic on $\mathring{C}_{n,q}$.The  point 1) of the theorem is proved.

 Let   $\mathring{f}$ be  $f$ restricted  to  $ \mathring{W}$. Then $ ( \pi \circ \nu \circ \mathring{f}): ( \mathring{W},p) \to (\C ^2,0) $ is a finite morphism which is analytically isomorphic to  $\pi \circ \nu $. But  the matrix associated to $\pi \circ \nu $ is the matrix $A_1$ computed in \ref{equotient}. The morphism  $\psi=(\pi \circ \nu \circ \mathring{f}) $ has generic order $n$ and has degree $1$ on its  ramification locus. The existence of the analytic isomorphism $\mathring{f}$ implies the unicity of $\psi$. The point 2) is proved.

 {\it{  end of proof of \ref{th1}}}

\noindent {\bf Adresses.}

\vskip.1in

\noindent Fran\c coise Michel / Laboratoire de Math\' ematiques  / IMT Toulouse /
Universit\' e Paul Sabatier / 118 route de Narbonne / F-31062
Toulouse / FRANCE

e-mail: fmichel@math.univ-toulouse.fr

\noindent Caude Weber / Institut de  Math\' ematiques /Universit\'e de Gen\`eve / Gen\`eve  / Suisse

e-mail: webernajac@gmail.com

\printindex

\end{document}